\rmfamily%

\documentclass [10pt]{article}
\usepackage{amsmath}
\usepackage{amsfonts}
\usepackage{amssymb}
\usepackage{graphicx}
\usepackage{float}

\newtheorem{lemma}{Lemma}
\newtheorem{conjecture}{Conjecture}

\newtheorem{corollary}{Corollary}

\begin{document}
	
\title{The Farey Sequence and the Mertens Function}
	
\author{Darrell Cox, Sourangshu Ghosh, and Eldar Sultanow}
	
\date{May 22, 2021}
	
\maketitle
	
\begin{abstract}	

\noindent Franel and Landau derived an arithmetic statement involving the Farey sequence that is equivalent to the Riemann hypothesis. Since there is a relationship between the Mertens function and the Riemann hypothesis, there should be a relationship between the Mertens function and the Farey sequence. Functions of subsets of the fractions in Farey sequences that are analogous to the Mertens function are introduced. Mikol\'{a}s proved that the sum of certain Mertens function values is 1. Results analogous to Mikol\'{a}s' theorem are the defining property of these functions. A relationship between the Farey sequence and the Riemann hypothesis other than the Franel-Landau theorem is postulated.  This conjecture involves a theorem of Mertens and the second Chebyshev function.

\end{abstract}

\section{Introduction}

Mikol\'{a}s~\cite{l} proved that $\sum_{n=1}^{x}M(\lfloor x/n \rfloor)=1$ where $M$ denotes the Mertens function. If $f$ and $g$ are two arithmetical function, their Dirichlet product (or convolution) is defined by the equation $h(n)=\sum_{d|n}f(d)g(n/d)$ (denoted by $f*g$ for $h$ and $(f*g)(n)$ for $h(n)$). Let $F$ denote a real or complex-valued function defined on the positive real axis (0, +$\infty$) such that $F(x)=0$ for $0<x<1$. A more general convolution is $\sum_{n\le x}\alpha(n)F(x/n)$ where $\alpha$ is any arithmetical function. The sum defines a new function $G$ on $(0, +\infty)$ which also vanishes for $0<x<1$. The function $G$ is denoted by $\alpha \circ F$. An associative property relating $\circ$ and $*$ is given by Lemma~\ref{la:1}.

\begin{lemma}
	\label{la:1}
	For any arithmetical function $\alpha$ and $\beta$, $\alpha \circ (\beta \circ F)=(\alpha * \beta)\circ F$.
\end{lemma}

\noindent (See section 2.14 of Apostol's~\cite{ta} book.) The functions considered here are $\alpha(n)=1$, $b(n)=\lfloor n/2\rfloor -2\lfloor n/4 \rfloor$ and $F(n)=M(n)$. The Farey sequence of order $n$ (denoted by $\mathcal{F}_{n})$ is the ascending sequence of irreducible fractions between 0 and 1 whose denominators do not exceed $n$. Corollary 1.33 in Matveev's~\cite{mv} book is given by Lemma~\ref{la:2}

\begin{lemma}
	\label{la:2}
	If $f_{t} \in \mathcal{F}_{n}-\{1/1\}$ where $n>1$, then $t=\sum_{j=2}^{n}M(n/j)\lfloor jf_{t}\rfloor$ where $M(\cdot)$ is the Mertens function.
\end{lemma}

\noindent Previously, Mikol\'{a}s proved a slightly more general version of this. Let $\varrho_{\nu}$ denote the $\nu$-th fraction in $F_{x}$ (the Farey sequence of order $x$) and $\mu(n)$ the M\"obius function. Mikol\'{a}s' Lemma 4 is given by Lemma~\ref{la:3}.

\begin{lemma}
	\label{la:3}
	Let $0\le\xi\le 1$ and let us denote by $h(\xi,x)$ the number of fractions in $F_{x}$ which are not greater than $\xi$. Then we have $h(\xi,x)=\sum_{\varrho_{\nu}\le\xi}1=\sum_{n=1}^{[x]}[n\xi]M(\frac{x}{n})=\sum_{n=1}^{[x]}\mu(n)\sum_{d=1}^{[x/n]}[d\xi]$.
\end{lemma}

\noindent $\beta \circ M$ is then the difference between the number of fractions less than or equal to 1/2 and greater than 1/4 and the number of fractions less than or equal to 1/4. 

\section{The Farey Sequence and Redheffer matrices}

Mikol\'{a}s proved that $\sum_{i=1}^{x}M(\lfloor x/i \rfloor)=1$. In general, $\sum_{i=1}^{x}M(\lfloor x/(in) \rfloor)=1$, $n=1,2,3,\ldots,x$ (since $\lfloor\lfloor x/n \rfloor/i\rfloor=\lfloor x/(in) \rfloor$). Let $R'$ denote a square matrix where element $(i, j)$ equals 1 if $j$ divides $i$ or 0 otherwise. In a Redheffer matrix, element $(i, j)$ equals 1 if $i$ divides $j$ or if $j=1$. Redheffer~\cite{r} proved that the determinant of such a $x$ by $x$ matrix equals $M(x)$. Let $T$ denote the matrix obtained from $R'$ by element-by-element multiplication of the columns by $M(\lfloor x/1 \rfloor),M(\lfloor x/2 \rfloor),M(\lfloor x/3 \rfloor),\ldots,M(\lfloor x/x \rfloor)$. For example, the $T$ matrix for $x=12$ is \\

\begin{tabular}{rrrrrrrrrrrr}
	$-2$ & 0 & 0 & 0 & 0 & 0 & 0 & 0 & 0 & 0 & 0 & 0 \\
	$-1$ & $-1$ & 0 & 0 & 0 & 0 & 0 & 0 & 0 & 0 & 0 & 0 \\
	$-1$ & 0 & $-1$ & 0 & 0 & 0 & 0 & 0 & 0 & 0 & 0 & 0 \\
	$-1$ & $-1$ & 0 & $-1$ & 0 & 0 & 0 & 0 & 0 & 0 & 0 & 0 \\
	0 & 0 & 0 & 0 & 0 & 0 & 0 & 0 & 0 & 0 & 0 & 0 \\
	0 & 0 & 0 & 0 & 0 & 0 & 0 & 0 & 0 & 0 & 0 & 0 \\
	1 & 0 & 0 & 0 & 0 & 0 & 1 & 0 & 0 & 0 & 0 & 0 \\
	1 & 1 & 0 & 1 & 0 & 0 & 0 & 1 & 0 & 0 & 0 & 0 \\
	1 & 0 & 1 & 0 & 0 & 0 & 0 & 0 & 1 & 0 & 0 & 0 \\
	1 & 1 & 0 & 0 & 1 & 0 & 0 & 0 & 0 & 1 & 0 & 0 \\
	1 & 0 & 0 & 0 & 0 & 0 & 0 & 0 & 0 & 0 & 1 & 0 \\
	1 & 1 & 1 & 1 & 0 & 1 & 0 & 0 & 0 & 0 & 0 & 1 
\end{tabular} \\

\noindent Let $\Phi(x)$ denote $\sum_{i=1}^{x}\varphi(i)$ where $\varphi$ is Euler's totient function. This is also the number of fractions in a Farey sequence of order $x$. Let $U$ denote the matrix obtained from $T$ by element-by-element multiplication of the columns by $\varphi(j)$. The sum of the sums of the columns of $U$ then equals $\Phi(x)$. $i=\sum_{d|i}\varphi(d)$, so $\sum_{i=1}^{x}M(\lfloor x/i \rfloor)i$ (the sum of the sums of the rows of $U$) equals $\Phi(x)$. This relationship given by equation~\ref{eq:4} was proved by Cox~\cite{c}. 

\begin{equation}
	\label{eq:4}
	\sum_{i=1}^{x}M(\lfloor x/i \rfloor)i=\Phi(x)
\end{equation}

\noindent Previously, Mikol\'{a}s proved this using a different approach. Let $f(n)$ and $g(n)$ be arbitrary arithmetical functions. Mikol\'{a}s Lemma 2 is given by equation~\ref{eq:5}.

\begin{equation}
	\label{eq:5}
	\sum_{n=1}^{[x]}\sum_{d|n}f(d)g(\frac{n}{d})=\sum_{d=1}^{[x]}f(d)\sum_{\delta=1}^{[x/d]}g(\delta)=\sum_{d=1}^{[x]}g(d)\sum_{\delta=1}^{[x/d]}f(\delta)
\end{equation}

\noindent Special cases of this are given by equation~\ref{eq:6}.

\begin{equation}
	\label{eq:6}
	\sum_{n=1}^{[x]}\sum_{d|n}f(d)=\sum_{d=1}^{[x]}[\frac{x}{d}]f(d)=\sum_{d=1}^{[x]}\sum_{\delta=1}^{[x/d]}f(\delta)
\end{equation}

\noindent From this, Mikol\'{a}s determined that equations~\ref{eq:7}~and~\ref{eq:8} hold.

\begin{equation}
	\label{eq:7}
	\sum_{n=1}^{[x]}M(\frac{x}{n})=\sum_{n=1}^{[x]}\mu(n)[\frac{x}{n}]=1
\end{equation}

\begin{equation}
	\label{eq:8}
	\Phi(x)=\sum_{n=1}^{[x]}\varphi(n)=\sum_{n=1}^{[x]}nM(\frac{x}{n})=\frac{1}{2}\sum_{n=1}^{[x]}\mu(n)[\frac{x}{n}]^{2}+\frac{1}{2}
\end{equation}

\noindent Let $\sigma_{x}(i)$ denote the sum of positive divisors function ($\sigma_{x}(i)=\sum_{d|i}d^{x}$). Another relationship proved by Cox is given by Lemma~\ref{la:4}.

\begin{lemma}
	\label{la:4}
	$\sum_{i=1}^{x}M(\lfloor x/i \rfloor)\sigma_{0}(i)=x$
\end{lemma}

\noindent The right-hand side of the equation is simply the number of columns in the above modified Redheffer matrix. The sum of each column is 1 and the respective sums of the rows are $M(\lfloor x/1 \rfloor)\sigma_{0}(1)$, $M(\lfloor x/2 \rfloor)\sigma_{0}(2)$, $M(\lfloor x/3 \rfloor)\sigma_{0}(3)$,$\ldots$,$M(\lfloor x/x \rfloor)\sigma_{0}(x)$.

\section{Properties of $1 * \beta$}

\noindent  A result from the literature (see 2.11 of Tenenbaum's~\cite{te} book) is given by Lemma~\ref{la:5}.
\begin{lemma} 
	\label{la:5} 
 $1*1=\sigma_{0}$ 
\end{lemma}

\noindent A more general result from the literature (see 2.7 of Tenenbaum's book) is given by Lemma~\ref{la:6}.
\begin{lemma}
	\label{la:6}
\noindent $(1*f)(n)=\sum_{d|n}f(d)$ 
\end{lemma}

\begin{figure}[H]
	\includegraphics[clip, trim=0cm 1cm 0cm 0cm, width=\linewidth]{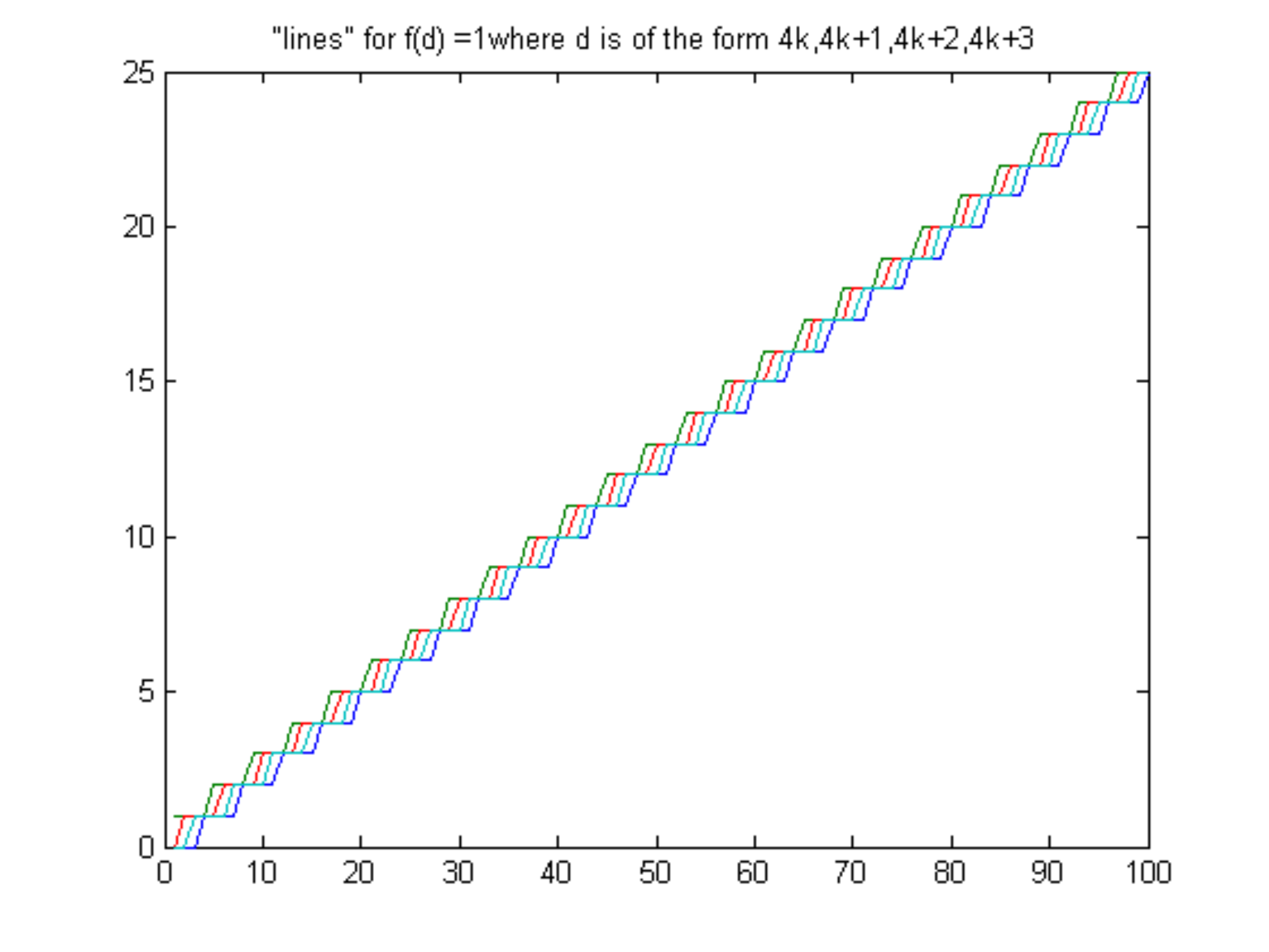}
	\caption{Plot of lines resulting from convolution with Mertens function}
	\label{fig:1}
\end{figure}

\noindent Note that $\lfloor j/2 \rfloor - 2\lfloor j/4 \rfloor$ equals $\lfloor (j+4)/2 \rfloor -2 \lfloor (j+4)/4 \rfloor$. The values of $\lfloor j/2 \rfloor - 2\lfloor j/4 \rfloor$ for $j=1,2,3,\ldots$ are then $0,1,1,0,0,1,1,0,0,1,1,0\ldots$. The values depend on whether $j$ is of the form $4k$, $4k+1$, $4k+2$, or $4k+3$. Let $f_{1}(d)=1$ if $d=4k$, $f_{2}(d)=1$ if $d=4k+1$, $f_{3}(d)=1$ if $d=4k+2$, and $f_{4}(d)=1$ if $d=4k+3$. Then by Lemma~\ref{la:4}, Lemma~\ref{la:5}, and Lemma~\ref{la:6}, there are four ``lines" with a slope of $\frac{1}{4}$ (out of phase by one position) corresponding to the convolution of these functions with the Mertens function. A plot of these ``lines" is given in Figure~\ref{fig:1}.

\noindent Note that two of the above ``lines" can be added to give a ``line" with a slope of $\frac{1}{2}$, three can be added to give a ``line" with a slope of $\frac{3}{4}$, or all of them can be added to give a straight line with a slope of 1. \\

\noindent Let $m_{x}$ denote the number of fractions up to and including $\frac{1}{4}$ and $n_{x}$ the number of fractions after $\frac{1}{4}$ and up to and including $\frac{1}{2}$ in a Farey sequence of order $x$ ($m_{1}$, $m_{2}$, $n_{1}$ and $n_{2}$ are set to 0). A consequence of the above convolution is given by Corollary~\ref{co:1}.

\begin{corollary}
	\label{co:1}
	$\sum_{i=1}^{x}(n_{\lfloor x/i \rfloor}-m_{\lfloor x/i \rfloor}+\frac{1}{4})$ equals $\frac{1}{2}$, $\frac{1}{4}$, $0$, or $-\frac{1}{4}$
\end{corollary}

\section{Farey Sequences and the Riemann Hypothesis}

For $\nu=1,2,3,\ldots,\Phi(x)$ let $\delta_{\nu}$ denote the amount by which the $\nu$-th term of the Farey sequence differs from $\nu/\Phi(x)$. The theorem of Franel and Landau~\cite{fl} is that the Riemann hypothesis is equivalent to the statement that $|\delta_{1}|+|\delta_{2}|+|\delta_{3}|+,\ldots,+|\delta_{\Phi(x)}|=o(x^{1/2+\epsilon})$ for all $\epsilon>0$ as $x\rightarrow \infty$. \\

\noindent Let $\Lambda(i)$ denote the Mangoldt function ($\Lambda(i)$ equals $\log(p)$ if $i=p^{m}$ for some prime $p$ and some $m\ge 1$ or $0$ otherwise). Let $\psi(x)$ denote the second Chebyshev function ($\psi(x)=\sum_{i\le x}\Lambda(i)$). Mertens~\cite{m1} proved that equation~\ref{eq:14} holds.

\begin{equation}
\label{eq:14}
\sum_{i=1}^{x}M(\lfloor x/i \rfloor)\log(i)=\psi(x)
\end{equation}

\begin{figure}[H]
	\includegraphics[clip, trim=0cm 1cm 0cm 0cm, width=\linewidth]{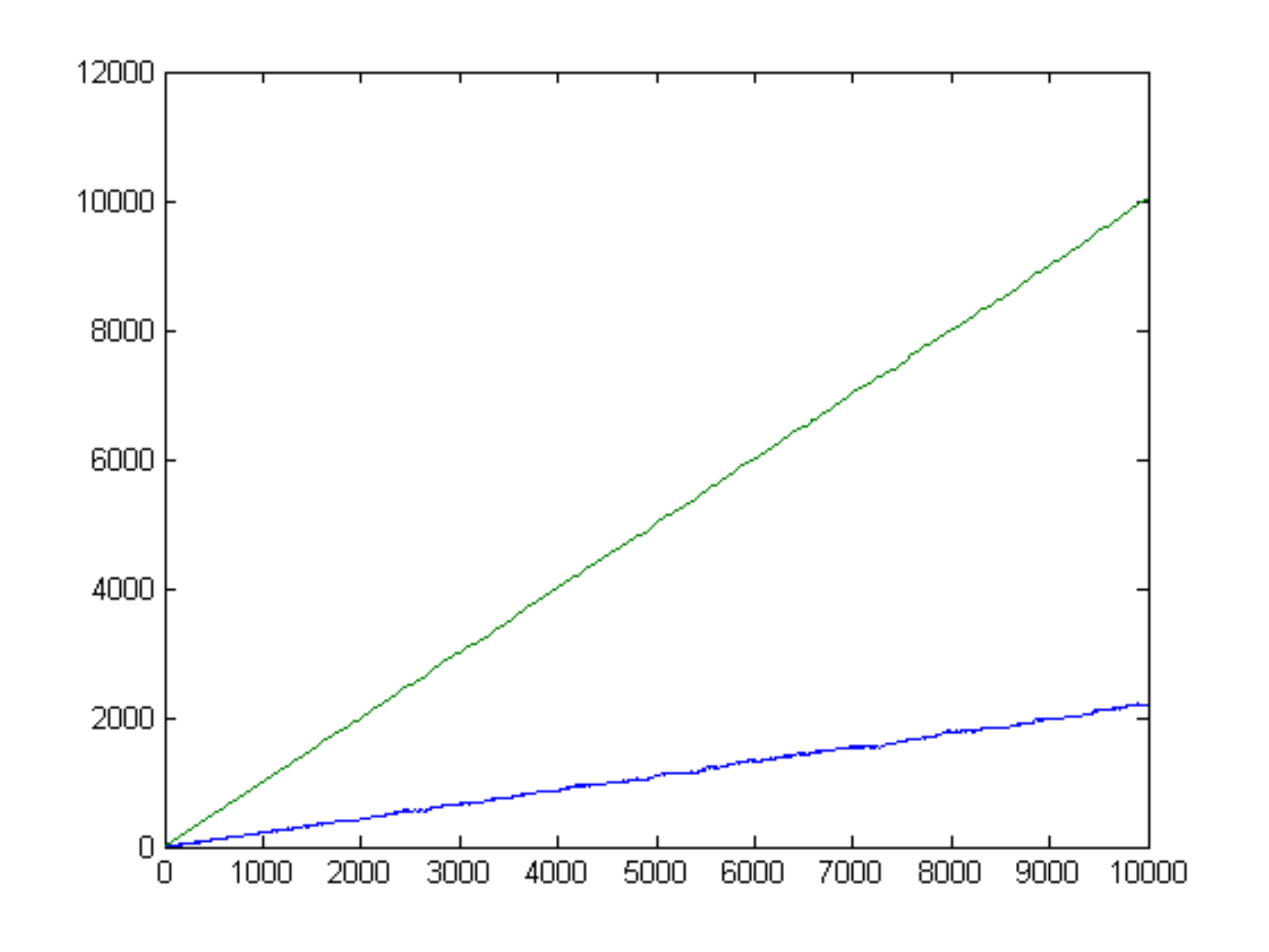}
	\caption{Plot of the second Chebyshev function and corresponding sums }
	\label{fig:2}
\end{figure}

\noindent A plot of $\psi(x)$ and $\sum_{i=1}^{x}(m_{\lfloor x/i \rfloor}-n_{\lfloor x/i\rfloor}+\frac{1}{4})\log(i)$ for $x=2,3,4,\ldots,10000$ is given in Figure~\ref{fig:2}.

\noindent The slope of $\sum_{i=1}^{x}(m_{\lfloor x/i \rfloor}-n_{\lfloor x/i\rfloor}+\frac{1}{4})\log(i)$ is approximately 0.2197 whereas the slope of $\psi(x)$ is approximately 1.0. \\

\noindent The Riemann hypothesis is equivalent to the arithmetic statement $\psi(x)-x=o(x^{1/2+\epsilon})$ for all $\epsilon>0$. For a linear least-squares fit of $\psi(x)$ for $x=2,3,4,\ldots,10000$, $p_{1}=1$ with a 95\% confidence interval of (1, 1), $p_{2}=-2.122$ with a 95\% confidence interval of ($-2.699$, $-1.545$), SSE=$2.165\cdot10^{6}$, R-squared=1, and RMSE=14.72. For a linear least-squares fit of 
$\sum_{i=1}^{x}(m_{\lfloor x/i \rfloor}-n_{\lfloor x/i\rfloor}+\frac{1}{4})\log(i)$ for $x=2,3,4,\ldots,10000$, $p_{1}=0.2197$ with a 95\% confidence interval of (0.2196, 0.2198), $p_{2}=-1.246$ with a 95\% confidence interval of ($-1.841$, $-0.6504$), SSE=$2.303\cdot10^6$, R-squared=0.9994, and RMSE=15.18. For a linear least-squares fit of $\psi(x)$ for $x=2,3,4,\ldots,25000$, $p_{1}=1$ with a 95\% confidence interval of (1, 1), $p_{2}=-5.18$ with a 95\% confidence interval of ($-5.821$, $-4.54$), SSE=$1.667\cdot10^{7}$, R-squared=1, and RMSE=25.82. For a linear least-squares fit of 
$\sum_{i=1}^{x}(m_{\lfloor x/i \rfloor}-n_{\lfloor x/i\rfloor}+\frac{1}{4})\log(i)$ for $x=2,3,4,\ldots,25000$, $p_{1}=0.2197$ with a 95\% confidence interval of (0.2197, 0.2198), $p_{2}=-2.773$ with a 95\% confidence interval of ($-3.358$, $-2.189$), SSE=$1.388\cdot10^7$, R-squared=0.9998, and RMSE=23.56. The sum-squared errors (SSE) and root-mean-squared errors (RMSE) are approximately equal. This is the basis for the following Conjecture~\ref{cj:1}.

\begin{conjecture}
\label{cj:1}
The Riemann hypothesis is equivalent to the arithmetic statement $\sum_{i=1}^{x}(m_{\lfloor x/i \rfloor}-n_{\lfloor x/i\rfloor}+\frac{1}{4})\log(i)-0.2197x=o(x^{1/2+\epsilon})$ for all $\epsilon>0$.
\end{conjecture}

\section{Infinitely Many Analogues of the Mertens Function}

Let $o_{x}$ denote the number of fractions up to and including $\frac{1}{5}$ and $p_{x}$ the number of fractions after $\frac{1}{5}$ and up to and including $\frac{2}{5}$ in a Farey sequence of order $x$. Let $q_{x}$ denote the number of fractions up to including $\frac{1}{6}$ and $r_{x}$ the number of fractions after $\frac{1}{6}$ and up to and including $\frac{1}{3}$ in a Farey sequence of order $x$. Let $s_{x}$ denote the number of fractions up to and including $\frac{1}{7}$ and $t_{x}$ the number of fractions after $\frac{1}{7}$ and up to and including $\frac{2}{7}$ in a Farey sequence of order $x$. Let $u_{x}$ denote the number of fractions up to and including $\frac{1}{8}$ and $v_{x}$ the number of fractions after $\frac{1}{8}$ and up to and including $\frac{1}{4}$ in a Farey sequence of order $x$. The following corollaries can be proved by using Lemma~\ref{la:1}, Lemma~\ref{la:2}, Lemma~\ref{la:4}, Lemma~\ref{la:5}, and Lemma~\ref{la:6} just as they were used to prove Corollary~\ref{co:1}.   

\begin{corollary}
	\label{co:2}
	$\sum_{i=1}^{x}(o_{\lfloor x/i \rfloor}-p_{\lfloor x/i \rfloor}+\frac{2}{5})$ equals $\frac{4}{5}$, $\frac{1}{5}$, $0$, $-\frac{1}{5}$, or $-\frac{2}{5}$
\end{corollary}

\noindent For a linear least-squares fit of $\sum_{i=1}^{x}(o_{\lfloor x/i \rfloor}-p_{\lfloor x/i \rfloor}+\frac{2}{5})\log(i)$ for $x=2,3,4,\ldots,25000$, $p_{1}=0.3733$ with a 95\% confidence interval of (0.3733, 0.3734), $p_{2}=1.598$ with a 95\% confidence interval of (0.9739, 2.222), SSE=$1.583\cdot10^{7}$, R-squared=0.9999, and RMSE=25.17. 

\begin{corollary}
	\label{co:3}
	$\sum_{i=1}^{x}(q_{\lfloor x/i \rfloor}-r_{\lfloor x/i \rfloor}+\frac{1}{3})$ equals $1$, $\frac{2}{3}$, $\frac{1}{3}$, $0$, or $-\frac{1}{3}$
\end{corollary}

\noindent For a linear least-squares fit of $\sum_{i=1}^{x}(q_{\lfloor x/i \rfloor}-r_{\lfloor x/i \rfloor}+\frac{1}{3})\log(i)$ for $x=2,3,4,\ldots,25000$, $p_{1}=0.4218$ with a 95\% confidence interval of (0.4218, 0.4218), $p_{2}=-2.837$ with a 95\% confidence interval of ($-3.449$, $-2.225$), SSE=$1.523\cdot10^{7}$, R-squared=0.9999, and RMSE=24.69. 

\begin{corollary}
	\label{co:4}
	$\sum_{i=1}^{x}(s_{\lfloor x/i \rfloor}-t_{\lfloor x/i \rfloor}+\frac{3}{7})$ equals $\frac{9}{7}$, $\frac{6}{7}$, $\frac{5}{7}$, $\frac{3}{7}$, $\frac{1}{7}$, $0$, or $-\frac{3}{7}$
\end{corollary}

\noindent For a linear least-squares fit of $\sum_{i=1}^{x}(s_{\lfloor x/i \rfloor}-t_{\lfloor x/i \rfloor}+\frac{3}{7})\log(i)$ for $x=2,3,4,\ldots,25000$, $p_{1}=0.5572$ with a 95\% confidence interval of (0.5572, 0.5573), $p_{2}=-1.185$ with a 95\% confidence interval of ($-1.87$, $-0.5011$), SSE=$1.904\cdot10^{7}$, R-squared=1, and RMSE=27.6. 

\begin{corollary}
	\label{co:5}
	$\sum_{i=1}^{x}(u_{\lfloor x/i \rfloor}-v_{\lfloor x/i \rfloor}+\frac{3}{8})$ equals $\frac{9}{8}$, $\frac{3}{2}$, $\frac{7}{8}$, $\frac{3}{4}$, $\frac{3}{8}$, $\frac{1}{4}$, $0$, or $-\frac{3}{8}$
\end{corollary}

\noindent For a linear least-squares fit of $\sum_{i=1}^{x}(u_{\lfloor x/i \rfloor}-v_{\lfloor x/i \rfloor}+\frac{3}{8})\log(i)$ for $x=2,3,4,\ldots,25000$, $p_{1}=0.5782$ with a 95\% confidence interval of (0.5782, 0.5783), $p_{2}=0.9218$ with a 95\% confidence interval of ($0.2422$, $1.601$), SSE=$1.878\cdot10^{7}$, R-squared=1, and RMSE=27.41. \\

\noindent In general, the differences between the number of fractions up to and including $1/I$ and after $1/I$ and up to and including $2/I$, $I=4,5,6,\ldots$ are incremented by $\frac{\lfloor(I+1)/2 \rfloor -1}{I}$.  For $I=4$, $\lfloor j/2 \rfloor - 2\lfloor j/4 \rfloor$ equals $\lfloor (j+4)/2 \rfloor -2 \lfloor (j+4)/4 \rfloor$.  For $I>4$, similar expressions determine the periods and thus the slopes of the ``lines" and the increments required to compensate for the slopes. \\

\begin{figure}[H]
	\includegraphics[clip, trim=0cm 1cm 0cm 0cm, width=\linewidth]{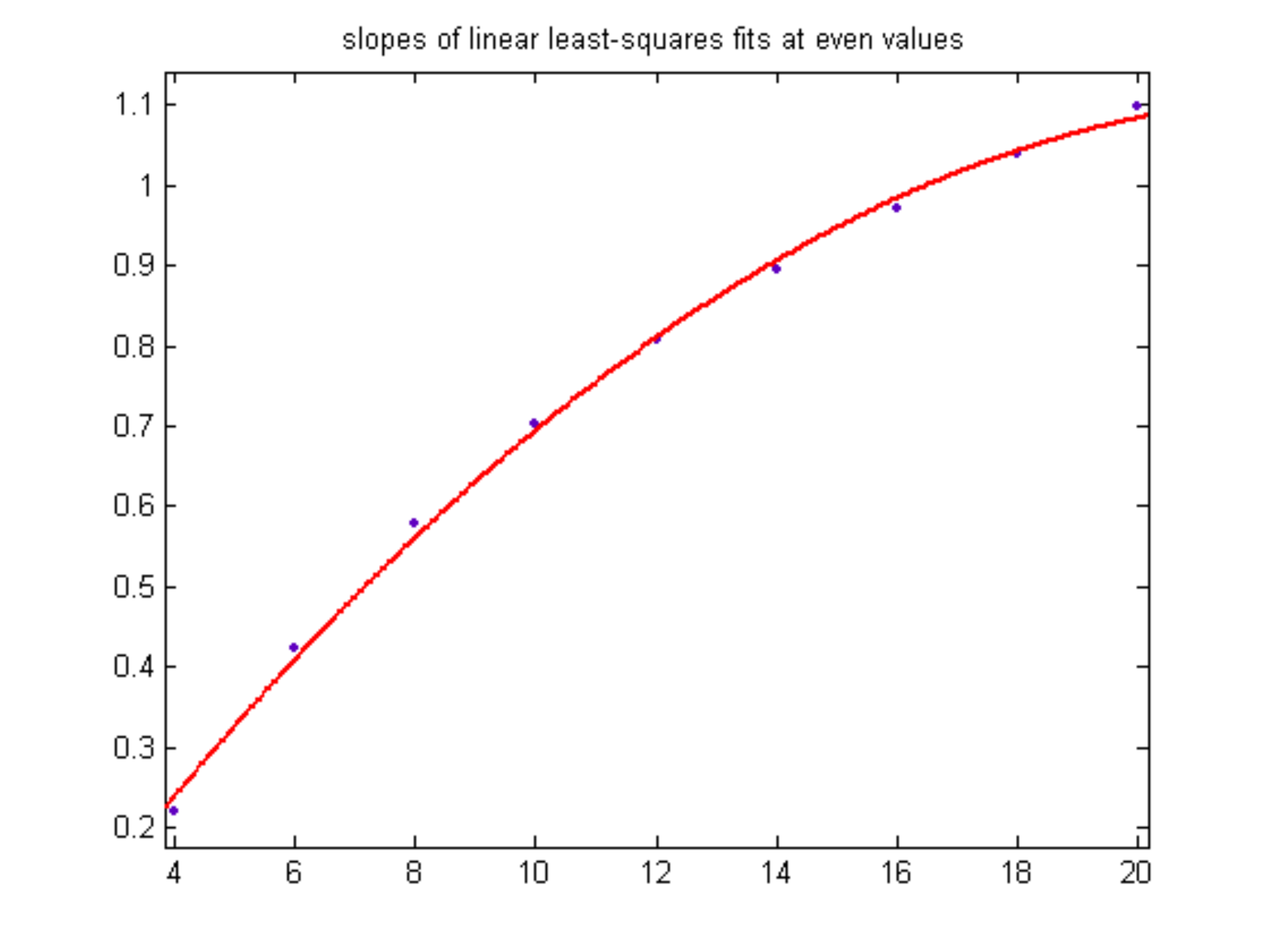}
	\caption{Plot of the slopes of linear least-squares fits at even I values}
	\label{fig:3}
\end{figure}

\noindent A plot of these and subsequent slopes at even $I$ values is given in Figure~\ref{fig:3}.

\begin{figure}[H]
	\includegraphics[clip, trim=0cm 1cm 0cm 0cm, width=\linewidth]{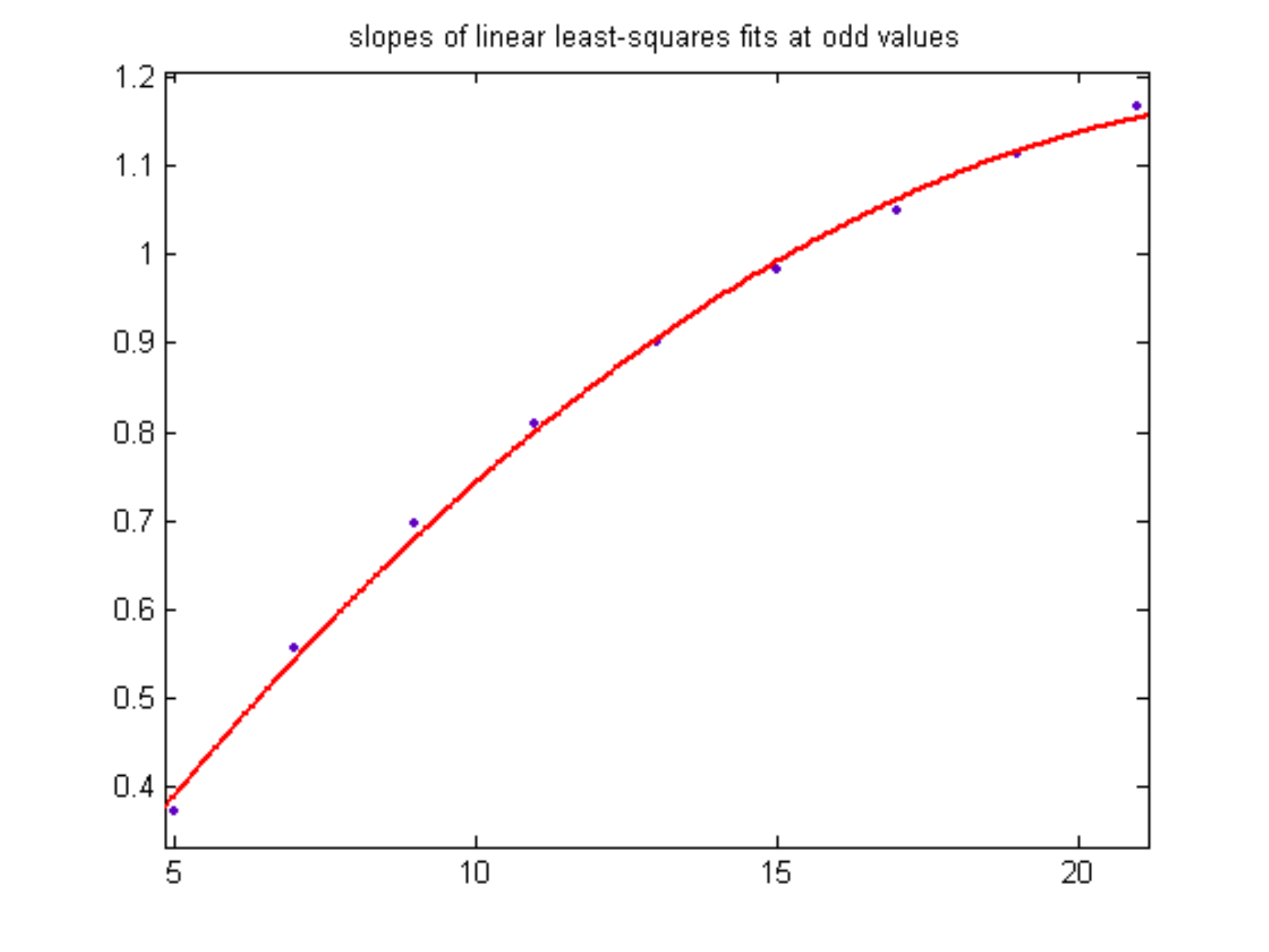}
	\caption{Plot of the slopes of linear least-squares fits at odd I values}
	\label{fig:4}
\end{figure}

\noindent For a quadratic least-squares fit of the data, R-squared=0.998. A plot of the slopes at corresponding odd $I$ values is given in Figure~\ref{fig:4}.

\noindent For a quadratic least-squares fit of the data, R-squared=0.9978.

\section{Conclusion}

\noindent The above analogues of the Mertens function can be substituted for the Mertens function in many equations and the results analyzed.  No method for proving the conjecture has been determined.

 \end{document}